\documentclass[reqno]{amsart}
\usepackage{amssymb}
\usepackage{upgreek}
\usepackage{xcolor}
\usepackage{tikz}

\theoremstyle{definition}

\numberwithin{equation}{section}
\allowdisplaybreaks
\begin{document}
\title[A century problem related to the Legendre symbol modulo $p$]
{A century problem related to the Legendre symbol modulo $p$}

\author{Zhang Wenpeng}
\address{Mathematics Research Institute,  Aba Teachers College, Wenchuan, Sichuan, China}
\curraddr{}
\email{wpzhang@nwu.edu.cn}


\date{October 16, 2026}

\maketitle

\begin{abstract}
The main purpose of this paper is using a very simple constructive method to study an old number theory problem related to the Legendre symbol modulo $p$, and completely solved it. The proving method of the result is purely elementary and has been desired in the literature at least since 1927.

\end{abstract}
\section*{1. Introduction}

Let $p$ be an odd prime. For any integer $a$, the famous Legendre symbol modulo $p$ is define (see [1] and [2]) as
\begin{eqnarray*}
\left(\frac{a}{p}\right)= \left\{
\begin{array}{ll} 1 &\textrm{ if $a$ is a quadratic residue modulo $p$;} \\   -1 &\textrm{ if $a$ is a quadratic non-residue modulo $p$;}\\
   0 &\textrm{ if $a\equiv 0$ modulo $p$.}
\end{array}\right.
\end{eqnarray*}
This number theory function occupies a very important position in the study of elementary number theory and analytic number theory, so many theorists and scholars have studied its various  properties, and a series of important results have been achieved. Perhaps the most significant result is its reciprocity formula. That is, for any two distinct odd primes $p$ and $q$, one has (see Theorem 9.8 in [1])
\begin{eqnarray*}
\left(\frac{q}{p}\right)\cdot \left(\frac{p}{q}\right)= (-1)^{\frac{(p-1)(q-1)}{4}}.
 \end{eqnarray*}
Some other properties of Legendre symbol can be found in references [3]--[6], and need not be listed here.

It is clear that
$$
\sum_{a=1}^{p-1}\left(\frac{a}{p}\right)=0.
$$
So in the reduced residue system modulo $p$, the number of the quadratic residues modulo $p$ and the number of the quadratic non-residues modulo $p$ half of each. If $p\equiv 1\bmod 4$, then note that $\left(\frac{-1}{p}\right)=1$, we have
$$
\sum_{a=1}^{p-1}\left(\frac{a}{p}\right)=\sum_{a=1}^{\frac{p-1}{2}}\left(\frac{a}{p}\right)+ \sum_{a=1}^{\frac{p-1}{2}}\left(\frac{p-a}{p}\right)=\sum_{a=1}^{\frac{p-1}{2}}\left(\frac{a}{p}\right)+\sum_{a=1}^{\frac{p-1}{2}}\left(\frac{-a}{p}\right)=2\sum_{a=1}^{\frac{p-1}{2}}\left(\frac{a}{p}\right).
$$
So we have the identity
\begin{eqnarray*}
A(p)=\sum_{a=1}^{\frac{p-1}{2}}\left(\frac{a}{p}\right)=0.
\end{eqnarray*}
That is to say, in the set $\{1,2,\cdots,\frac{p-1}{2}\}$, the number of the quadratic residues modulo $p$ and the number of the quadratic non-residues modulo $p$ also half of each.
But if $p\equiv 3\bmod 4$, then $\frac{p-1}{2}$ is an odd number and
\begin{eqnarray}
A(p)=\sum_{a=1}^{\frac{p-1}{2}}\left(\frac{a}{p}\right)\neq 0.
\end{eqnarray}

About a century ago,  J. Jacobi first conjecture that $A(p)>0$ and proved by L. Dirichlet [7] and [8] in connection with the theory of binary quadratic forms. The proofs are also given in the books of P. Bachmann [9] and E. Landau [10]. In particular, E. Landau and H. Davenport emphasizes in their books [10] and [15] that they hope to see a elementary proof for $A(p)>0$. More proofs are due to Kai-Lai Chung [11], A. L. Whiteman [12] and L. Moser [13]. It is a shame that all known proofs are analytic. So, L. Moser in [13] said with emotion, ``While a really elementary proof would be of great interest".

Through analytic methods, i.e. Dirichlet's class-number formula $h(p)=\frac{\sqrt{p}}{\pi}\cdot L(1,\chi_2)$, where $\chi_2=\left(\frac{*}{p}\right)$, $h(p)$ denotes the class number of the imaginary quadratic field $\mathbf{Q}\left(\sqrt{-p}\right)$ (see [2], p. 295 or [14], p. 395 or [15], p. 50), one can easy to prove that $A(p)>0$. In fact, from B. C. Berndt's work [16] (see Theorem 3.2) we have
\begin{eqnarray}
\sum_{a=1}^{\frac{p-1}{2}}\left(\frac{a}{p}\right)= \frac{i\cdot \tau\left(\chi_2\right)}{\pi}\cdot \left(\left(\frac{2}{p}\right)-2\right)\cdot L\left(1,\chi_2\right)=\frac{\sqrt{p}}{\pi}\left(2-\left(\frac{2}{p}\right)\right)\cdot L(1, \chi_2),
\end{eqnarray}
where $\tau(\chi_2)=i\cdot \sqrt{p}$.

 Then combining Dirichlet's class-number formula and (2) we have
 \begin{eqnarray*}
\sum_{a=1}^{\frac{p-1}{2}}\left(\frac{a}{p}\right)=\left(2-\left(\frac{2}{p}\right)\right)\cdot h(p).
\end{eqnarray*}
Since $2-\left(\frac{2}{p}\right)>0$ and $h(p)> 0$, so $A(p)>0$.

That is, if $p\equiv 3\bmod 4$, then in the set $\{1, \ 2, \ \cdots, \ \frac{p-1}{2}\}$, the number of the quadratic residues modulo $p$ is greater than the number of the quadratic non-residues modulo $p$. R. K. Guy [17] (see problem F5) asked us:

Is there an elementary proof for $A(p)>0$?

 About this problem, it seems to have not solved at present, at least we have not found any elementary proof in the existing literature. In theory, this problem can be proven by using the elementary methods, but no one has been able to find one. Of course, it would be very meaningful to have a simple elementary method to prove this problem, just like the proof of the prime number theorem, the problem was already solved over a century ago by using the non-zero region of the Riemann zeta-function. It is precisely because A. Selberg [18] proved the prime number theorem using elementary methods that he was awarded the Fields Medal. Although A. Wiles [19] proved Fermat's Last Theorem by using profound modular elliptic curves theory, I think it would be a remarkable achievement if someone could prove Fermat's Last Theorem  by using elementary methods. This is the way to make people feel that the problem and the proof method are aligned.

 The main purpose of this paper is using a very simple elementary constructive method and the properties of the Legendre symbol to prove the following:

{\bf Theorem.} For any odd prime $p$ with $p\equiv 3\bmod 4$, we have
$$
A(p)=\sum_{a=1}^{\frac{p-1}{2}}\left(\frac{a}{p}\right)>0.
$$

\section*{ 2. Proof of the theorem}

In this section, we provide a direct proof of the theorem. Of course, in the process of proving the theorem, we need to use some properties of the Legendre symbol, all them can be found in [1] or [2], and will not be repeated here.  Without loss of generality, we may suppose prime $p=4m+3> 31$. When $3< p< 11$, we can verify the theorem is correct one by one.  Let us distinguish two cases to discuss:

Case 1. $\left(\frac{2}{p}\right)=-1$. In this case, note that  $\left(\frac{2}{p}
\right)=(-1)^{\frac{p^2-1}{8}}$ and $\left(\frac{-1}{p}\right)=-1$, we have $p=8k+3$. It is clear that for any odd number $r$, $\left(\frac{(p-r)/2}{p}\right)=\left(\frac{r}{p}\right)$.
Let  $A= \{1, \ 2, \ 3, \ \cdots, \ 4k+1\}$. Now for any positive integers $2n+1\in A$, then  $0\leq n\leq 2k$. If  $2n+1$ is a quadratic residue modulo $p$, then $\left(p-2n-1\right)/2$ is also a quadratic residue modulo $p$, and  $\left(p-2n-1\right)/2 \in A$; $4(2n+1)\in A$ is a quadratic residue modulo $p$, when $0\leq n\leq \frac{k}{2}-\frac{3}{8}$; $4(2n+1)-p\in A$ is also a quadratic residue modulo $p$, when $k-\frac{1}{8}< n\leq \frac{3k}{2}$. In this case, we also have $\{2n+1: 0\leq n\leq \frac{4k}{3}\}\bigcap \{ \frac{p-2n-1}{2}:  0\leq n< \frac{4k}{3}\}=\phi$.

If $2n+1$ is not a quadratic residue modulo $p$, then  $2(2n+1)\in A$ is a quadratic residue modulo $p$, when $1\leq n<k$; $8(2n+1)\in A$ is a quadratic residue modulo $p$, when $1\leq n \leq \frac{k}{4}-\frac{7}{16}$; $8(2n+1)-p\in A$ is a quadratic residue modulo $p$, when $\frac{k}{2}-\frac{5}{16}<n\leq \frac{3k}{4}-\frac{1}{4}$; $p-4(2n+1)\in A$ is also a quadratic residue modulo $p$, when $\frac{k}{2}-\frac{1}{4}\leq n< k-\frac{1}{8}$.  Note that $1$, $4$, $2k+1$, $3k+1$ and $4k+1$ are quadratic residue modulo $p$. It is clear that in the interval $[1, \frac{4k-7}{16}]$, each $k$ corresponds to at least two quadratic residue modulo $p$; In the interval $[\frac{k}{4}-\frac{7}{16}, \frac{k}{2}-\frac{5}{16}]$,  each $k$ corresponds to at least one quadratic residue modulo $p$; In the interval $[\frac{k}{2}-\frac{1}{4}, \frac{3k}{4}-\frac{1}{4}]$, each $k$ corresponds to at least  three quadratic residue modulo $p$; In the interval $[\frac{3k}{4}-\frac{1}{4}, k-\frac{1}{8}]$, each $k$ corresponds to at least two quadratic residue modulo $p$. So when $1\leq n\leq  2k$, there are at least
\begin{eqnarray*}&&5+2\left[\frac{4k-7}{16}\right]+ \left[\frac{k}{2}+\frac{2}{16}-\frac{k}{4}\right]+ 3\left[\frac{3k}{4}-\frac{k}{2}\right]+ 2\left[k- \frac{3k}{4}+\frac{1}{8}\right]\nonumber\\
&=&5+2\left[\frac{4k-7}{16}\right]+\left[\frac{4k+2}{16}\right]+3\left[\frac{k}{4}\right]+ 2\left[ \frac{2k+1}{8}\right]\geq 2k+1
\end{eqnarray*}
quadratic residue modulo $p$ in the set $A$.

Therefore, for any prime $p>11$ with $p=8k+3$, we have
\begin{eqnarray}
A(p)=\sum_{a=1}^{\frac{p-1}{2}}\left(\frac{a}{p}\right)> 0.
\end{eqnarray}

 Case 2. $\left(\frac{2}{p}\right)=1$. In this case, for any odd number $r$ with $(r, p)=1$, we have $p=8k+7$ and $\left(\frac{(p-r)/2}{p}\right)=-\left(\frac{r}{p}\right)$. Let $A=\{1, 2, \cdots,4k+3\}$. For any $4n+3\in A$, it is clear that one of $4n+3$ or $4k-2n+2$ is a quadratic residue modulo $p$ in $A$. And $4n+3\in A$ if and only if $0\leq n\leq k$. So in this case, there are at least $k+1$ such quadratic residues in the set $A$.

 If $4n+3\leq 2k+1$, then when $n=2h+1$, both $4n+3$ and $2(4n+3)$ or $4k-4h$, $2k-2h$ and $k-h$ are three quadratic residues modulo $p$ in $A$. In this case, $4(2h+1)+3\leq 2k+1$ or $0\leq h\leq \frac{k-3}{4}$.  Thus, in the case of $4n+3$, we can construct at least
  $$
  k+1 + \left[\frac{k-3}{4}\right]+1= k+\left[\frac{k+1}{4}\right]+1
  $$
different quadratic residues modulo $p$ in $A$.

If $4n+1\in A$, then $0\leq n\leq k$, and one of $4n+1$ or $ \frac{p-4n-1}{2}= 4k-2n+3$ is a quadratic residue modulo $p$. So in this case, there are at least $k+1$ different quadratic residues modulo $p$ in $A$.

Now let us consider the possible duplication with the previous cases. For any old number $n=2h+1$, it is clear that  $4k-2n+3=4k-4h+1=4(k-h)+1$ is not repeated with the previous one. If $n=2h$, then we have $0\leq h\leq \frac{k}{2}$ and $4(k-h)+3\geq 2k+3$. This time, $4k-2n+3=4(k-h)+3$. If $4(k-h)+3$  is a quadratic residue modulo $p$, then $p-2(8h+1)= 8k-16h+5= 4(2k-4h+1)+1$ is also a quadratic residue modulo $p$. That is to say, if $4(2h)+1$ is not a quadratic residue modulo $p$, then  $p-2(8h+1)= 8(k-2h)+5= 4(2k-4h+1)+1$ is a quadratic residue modulo $p$. If $1\leq 2k-4h+1\leq k$, then $\frac{k+1}{4}\leq h\leq \frac{k}{2}$ and the quadratic residue $4(2k-4h+1)+1\in A$. So for $m=2k-4h+1$, $4m+1$ is a quadratic residue modulo $p$ and $4k-2m+3$ is not a quadratic residue modulo $p$. Therefore, in the possible repeating elements $8h+1$ and $4(k-h)+3$, $h$ must satisfy $0\leq h\leq \frac{k+1}{4}$. But for $h=0$, $1=8\cdot 0+1$ is a quadratic residues modulo $p$. Thus, in the case $4n+1$, we can construct at least
  $$
k+1- \left[\frac{k+1}{4}\right]
  $$
  different quadratic residues in the set $A$ that are different from the previous cases.

Note that the quadratic residues constructed above are not repeated.  Therefore, in the set $\{1, 2, 3, \cdots, 4k+2, 4k+3\}$, the number of all different quadratic residues modulo $p$ are at least
 \begin{eqnarray}
k+\left[\frac{k+1}{4}\right]+1+ k+1-\left[\frac{k+1}{4}\right]= 2k+2.
\end{eqnarray}
So if $p=8k+7$, then from (4) we also have
 \begin{eqnarray}
A(p)=\sum_{a=1}^{\frac{p-1}{2}}\left(\frac{a}{p}\right)> 0.
\end{eqnarray}

It is clear that $p\equiv 3\bmod 4$ if and only if $p=8k+3$ or $p=8k+7$.

Combining (3) and (5) we complete the proof of our theorem.

\bigskip

{\bf Acknowledgements:} The author would like to express his sincere gratitude to Professor Gong Ke, Firdavs Rakhmonov and He Bo for their comments and suggestions.

 \end{document}